\documentclass[reqno,a4paper,12pt]{amsart}
\usepackage{amsmath,amsthm,amsfonts,amssymb}
\usepackage{verbatim, graphicx, ifthen}
\usepackage[T1]{fontenc}

\allowdisplaybreaks

\let\oldmarginpar\marginpar
\renewcommand\marginpar[1]{\-\oldmarginpar[\raggedleft\footnotesize #1]%
{\raggedright\footnotesize #1}}
\newtheorem*{theorema}{Theorem A}
\newtheorem*{theoremb}{Theorem B}
\newtheorem*{theoremc}{Theorem C}
\newtheorem*{theoremd}{Theorem D}

\newtheorem*{theorem1'}{Theorem 1'}
\newtheorem*{theorem2'}{Theorem 2'}
\newtheorem*{theorem1''}{Theorem 1''}
\newtheorem*{theorem2''}{Theorem 2''}

\newtheorem{coro}{Corollary}

\newtheorem{clame}{Claim}

\newtheorem{theorem}{Theorem}

\theoremstyle{definition}

\newcommand{\Z}{\mathbb{Z}}

\newcommand{\N}{\mathbb{N}}
\newcommand{\R}{\mathbb{R}}
\newcommand{\T}{\mathbb{T}}
\newcommand{\C}{\mathbb{C}}

\def\L{\Lambda}
\def\l{\lambda}
\def\T{\mathbb{T}}
\def\N{\mathbb{N}}
\def\Z{\mathbb{Z}}

\def\R{\mathbb{R}}
\def\C{\mathbb{C}}
\def\F{\mathcal{F}}

\def\1{\mathbf{1}}
\def\eps{\varepsilon}

\begin{document}

 \title[A few remarks on Sampling of signals with small spectrum]{A few remarks on Sampling of signals with small spectrum}

\author{Shahaf Nitzan, Alexander Olevskii and Alexander Ulanovskii}

\begin{abstract}
Given a set $S$ of small measure, we discuss  existence of
 sampling sequences for the
    Paley-Wiener space $PW_S$, which have both densities and
    sampling bounds close to the optimal ones.
\end{abstract}

\maketitle

\section{Introduction}

Let $S\subset [0,2\pi ]$ be a set  of positive measure. Denote by
$PW_S$ the space of all functions $f\in L^2(\R)$ whose Fourier
transform, $$\hat f(t):=\int_\R e^{-ixt}f(x)\,\frac{dx}{2\pi},$$
is supported by $S$, endowed with the $L^2-$norm.

A
sequence $\L\subset \Z$ is called a sampling sequence for $PW_S$
if there exists $B>0$ such that
\begin{equation}\label{frame}
B\|f\|^2\leq \sum_{\l\in\L} |f(\l)|^2
\:\:\:\:\:\:\:\:\:\:\:\:\forall f\in PW_S.
\end{equation}

 Landau's theorem \cite{La2} states that sampling sequences cannot
 be too sparse compared to the measure of $S$: If $\L$ is a sampling sequence for $PW_S$, then
\[
D^-{(\L)}:=\lim_{r\rightarrow
\infty}\frac{\min_{x\in\R}|\L\cap[x,x+r]|}{r}\geq |S|,
\]
where $| \cdot |$ denotes the normalized Lebesgue measure on $[0,2\pi]$.

In this note we discuss the existence of "good" sampling
sequences $\L\subset\Z$. 
The following observation explains what we mean by that: Assume
$S=[0,2\pi/m ]$, where $m\in\N$. For every $j\in\N$, the sequence
$\L=\{j+mn:n\in\Z\}$ is an excellent sampling sequence for
$PW_{S}$: It is uniformly distributed in $\Z$, its density is
equal to $|S|$ (the minimal density allowed, by Landau's theorem)
and
\begin{equation}
|S|\|f\|^2=\sum_{\l\in \L} |f(\l)|^2.
\end{equation}

For arbitrary sets $S$, we show the
existence of sampling sequences $\L$ with both densities and
bounds close to the ones described in this canonical example. When
$S$ has a "large" measure, the existence of such sequences follows
from  Bourgain and Tzafriri's Restricted Invertibility Theorem
\cite{BT} (see also \cite{V}). However, this result apparently is
not  applicable in the case when the spectrum $S$ has a "small" measure (see Remark 1). Our results are based instead on
a recent finite-dimensional theorem by Batson, Spielman and
Srivastava \cite{S1}.

     Observe that in \cite{OU} Olevskii and Ulanovskii constructed universal sampling sequences with
     density near optimal, which provide sampling for all compacts of a given measure. For another approach to this result see \cite{MM}.
     However, it was proved in \cite{OU} that such sampling sequences never have "universally good" bounds.

We distinguish two cases: when $S$ is a compact, and when it is a general (measurable)
set. For a compact set $S$ there exists a  good sampling
sequence which is uniformly distributed in $\Z$. To be more
precise, set
\[
D^+{(\L)}:=\lim_{r\rightarrow
\infty}\frac{\max_{x\in\R}|\L\cap[x,x+r]|}{r}.
\]
If $D^-{(\L)}=D^+{(\L)}$, one says that $\L$ has uniform density
$D(\L):=D^-{(\L)}=D^+{(\L)}$.

\begin{theorem}
Fix $d>0$. Given a compact set $S\subset [0,2\pi ]$, there exists
$\L\subset\Z$ with $D(\L)< (1+d)|S|$ such that
\[
C(d)|S|\|f\|^2\leq \sum |f(\l)|^2 \:\:\:\:\:\:\:\:\:\:\:\:\forall
f\in PW_S,
\]
where $C(d)$ is a positive constant depending only on $d$.
\end{theorem}

Moreover, in this case the sequence $\L$ can be chosen as a
finite union of arithmetic progressions.

For general sets $S$ the situation is more delicate. In this case,
we do not know whether sequences $\L$ of uniform density can
provide good sampling sequences. However, if the uniform density
is replaced by a weaker notion of density:
\[
D^{\sharp}(\L):=\lim_{n\rightarrow\infty}\frac{\L\cap[-n,n]}{2n},
\]
(provided the limit exists), then the result holds.
\begin{theorem}
Fix $d>0$. Given a set $S\subset [0,2\pi ]$ of positive measure,
there exists $\L\subset\Z$ with $D^{\sharp}(\L)< (1+d)|S|$ such
that
\[
C(d)|S|\|f\|^2\leq \sum |f(\l)|^2 \:\:\:\:\:\:\:\:\:\:\:\:\forall
f\in PW_S,
\]
where $C(d)$ is a positive constant depending only on $d$.
\end{theorem}

As mentioned above, these theorems are based on a result in
\cite{S1}, which studied the existence of small, "well
invertible", submatrices. This result is presented in Section 2.
Theorem~1 is then proved in Section 3. In Section 4 we use Theorem~1 and a Theorem of Ruzsa \cite{R} to obtain Theorem 2. In Section
5 we discuss the upper bounds and in Section 6 some other related
questions and open problems.

Set $\T:=[0,2\pi]$.
Throughout the rest of this paper we assume that $\Lambda\subset\Z$ and
$S\subset \T$. We denote by $|S|$ the normalized Lebesgue measure of a set $S$, and by $|J|$
 the number of elements in a finite set $J$. We let $l^n_2$ be the $n$ dimensional
space over $\C$ with norm
$\|\{a_r\}_{r=1}^n\|^2=\sum_{r=1}^n|a_r|^2$. For $J\subset
\{1,...,n\}$ we denote by $l_2(J)$ the $|J|$ dimensional subspace
of $l^n_2$, indexed by $J$.  Given a matrix $A$ of order $m\times
n$, and a subset $J\subseteq \{1,...,m\}$, we denote by $A(J)$ the
sub-matrix of $A$ with rows belonging to the index set $J$.

\section{Batson, Spielman and Srivastava's theorem}\label{section: bss thm}

\begin{theorema}\label{thm: bss thm}
Let $\{v_i\}_{i=1}^{m}$ be a system of vectors in $l_2^n$, $n\leq
m$, which satisfies
\begin{equation}
\label{finite frame cond} \sum_{i=1}^{m}|\langle w, v_i
\rangle|^2=\|w\|^2\:\:\:\:\:\:\:\:\:\:\:\:\forall w\in l_2^n.
\end{equation}
 Then for every $d>0$
there exist, a subset $J\subset \{1,...,m\}$ and positive weights
$\{s_i\}_{i\in J} $, such that $|J|\leq (1+d) n$ and
\[
A(d)\|w\|^2\leq \sum_{i\in J} s_i |\langle w, v_i\rangle|^2\leq
B(d)\|w\|^2\:\:\:\:\:\:\:\:\:\:\:\:\forall w\in l_2^n,
\]
where $A(d)$ and $B(d)$ are positive constants depending only on
$d$.
\end{theorema}

Observe, that this theorem is formulated in \cite{S1} for the real
spaces. However, the proof works also for the complex spaces.

Assume additionally that $\|v_i\|^2={n}/{m} $ for every $i=1,...,m$. Then,
by putting $w=v_j$, $j\in J$, we find
\[
s_j |\langle v_j, v_j\rangle|^2\leq\sum_{i\in J} s_i |\langle v_j,
v_i\rangle|^2\leq B(d)\|v_j\|^2,\:\
\]
which implies that $s_j\|v_j\|^2\leq B(d)$ or $s_j\leq B(d)m/n$
for every j. It follows that
\[
C(d)\frac{n}{m}\|w\|^2\leq \sum_{i\in J} |\langle w,
v_i\rangle|^2\:\:\:\:\:\:\:\:\:\:\:\:\forall w\in l_2^n
\]
where $C(d)$ is a positive constant depending only on $d$. Hence,
we get
\begin{coro}\label{coro: bss coro}
Let $A$ be an $m\times n$ matrix which is a sub-matrix of some
$m\times m$ orthonormal matrix, and such that all of it's rows
have equal $l^2$ norm. Then for every $d>0$ there exists a subset
$J\subset \{1,...,m\}$ for which $|J|\leq (1+d) n$ and
\[
C(d)\frac{n}{m}\|w\|^2\leq
\|A(J)w\|_{l_2(J)}^2\:\:\:\:\:\:\:\:\:\:\:\:\forall w\in l_2^n,
\]
where $C(d)$ is a positive constant depending only on $d$.
\end{coro}

\section{Proof of Theorem 1}\label{Section: proof of thm 1}

\noindent
Fix numbers $n,m\in\N$ and $d>0$ satisfying  $n(1+d)\leq m.$

\medskip\noindent
3.1. It suffices to prove Theorem 1 for the sets $S$ of the form
$$
S=\bigcup_{r\in I}\left[\frac{2\pi r}{m},\frac{2\pi(r+1)}{m}\right],
$$
where $I\subset\{0,...,m-1\}$, $|I|=n.$ Clearly, $|S|=n/m.$

\medskip\noindent
3.2. Denote by $$\F_I:=(e^{i\frac{2\pi jr}{m}})_{r\in
I,j=0,...,m-1 }$$ the submatrix of the Fourier matrix $\F$ whose
columns are indexed by $I$. Since the matrix $(\sqrt m)^{-1}\F$ is
orthonormal, by Corollary~1 there exists $J\subset\{0,...,m-1\}$,
$|J|\leq (1+d)n$, such that
\begin{equation}\label{f}
\Vert \F_I(J)w\Vert^2_{l_2(J)}\geq C(d)n\Vert w\Vert^2,
\:\:\:\:\:\:\:\:\:\:w\in l_2(I).\end{equation}

\noindent 3.3. Observe that every function $F\in L^2(S)$ can be
written as $$F(t)=\sum_{r\in I}F_r(t-\frac{2\pi r}{m}), $$ where
$F_r\in L^2(0,\frac{2\pi}{m})$ is defined by
$$F_r(t):=F(t+\frac{2\pi r}{m}){\bf 1}_{[0,\frac{2\pi}{m}]}(t).$$ Therefore, every function $f\in PW_S$ admits a representation
$$
f(x)=\sum_{r\in I}e^{i\frac{2\pi
r}{m}x}f_r(x),\:\:\:\:\:\:\:\:\:\: f_r\in PW_{[0,
\frac{2\pi}{m}]},
$$
where the functions $e^{i\frac{2\pi r}{m}x}f_r(x)$ are orthogonal
in $L^2(\R)$.

\noindent 3.4. We now verify that the sequence
$$
\L:=\{j+km: j\in J, k\in\Z\}
$$
satisfies the conclusion of Theorem 1. Take any function $f\in PW_S$. Then
$$
\sum_{j\in J}\sum_{k\in\Z}|f(j+km)|^2=\sum_{j\in
J}\sum_{k\in\Z}\left|\sum_{r\in I}e^{i\frac{2\pi rj
}{m}}f_r(j+km)\right|^2.
$$
  For every $j\in J$ we apply (2) to the function $\sum_{r\in I}e^{i\frac{2\pi rj
}{m}}f_r(x)$. We find that the last expression is equal to
$$
\frac{1 }{m}\sum_{j\in J}\int_{\R}\left|\sum_{r\in
I}e^{i\frac{2\pi rj }{m}}f_r(x)\right|^2dx= \frac{1
}{m}\int_{\R}\|\F_I(J)(f_r(x))_{r\in I}\|^2_{l_2(J)}dx.$$ By
inequality (4), we have
$$
\sum_{\l\in\L}|f(\l)|^2\geq C(d)\frac{n }{m}\int_{\R}\sum_{r\in
I}|f_r(x)|^2dx=
$$$$C(d)\frac{n }{m}\int_{\R}\sum_{r\in I}|e^{i\frac{2\pi
r}{m}x}f_r(x)|^2dx= C(d)\frac{n }{m}\int_{\R}|\sum_{r\in
I}e^{i\frac{2\pi r}{m}x}f_r(x)|^2dx=$$$$ C(d)\frac{n
}{m}\int_{\R}|f(x)|^2dx.
$$
This completes the proof.

\section{Proof of Theorem 2}
\subsection{Auxiliary results} We will use the following theorem of Ruzsa, \cite{R}
(see also \cite{V} for the extension to a two-sided density).

\begin{theoremb}
Let $\textit{H}$ be a family of finite sets of integers which
 is closed to translations, and
such that all subsets of a set in $\textit{H}$ also belong to
$\textit{H}$. Set
\[
d(\textit{H}):=\lim_{n\rightarrow\infty}\frac{\max_{A\in\textit{H}}|A\cap[-n,n]|}{2n}.
\]
Then there exists $\Gamma\subset\Z$ with
$D^{\sharp}(\Gamma)=d(\textit{H})$ such that every finite subset of
$\Gamma$ belongs to $\textit{H}$.
\end{theoremb}

We note that the limit $d(\textit{H})$ always exists.

Theorem B was used by Bourgain and Tzafriri in an application of
the Restricted Invertibility Theorem. Here we use it in a
different way,
           in order to manage the case of sampling
           with small spectrum.

As Theorem B studies families of finite sets, it is easier applied
when working with interpolating sequences (equivalently, Riesz sequences), rather than sampling sequences.
A system $\{e^{i\gamma t}\}_{\gamma\in \Gamma}$ is called a Riesz
sequence for $L^2(\Omega)$ if
\begin{equation} \label{rs}
A\sum_{\gamma\in\Gamma}|a_{\gamma}|^2\leq \|\sum_{\gamma\in
\Gamma} a_{\gamma}e^{i\gamma
t}\|_{L^2(\Omega)}^2\:\:\:\:\:\:\:\:\:\:\:\:\:\:\forall\{a_{\gamma}\}_{\gamma\in\Gamma}\in
l^2(\Gamma) ,
\end{equation}
where $A$ is a positive constant not depending on $a_{\gamma}$. The
duality between
 Riesz and sampling sequences is well known:

\begin{clame}A sequence $\L$ is a sampling sequence for $PW_S$ if and
only if $\{e^{i\gamma t}\}_{\gamma\in\Z\setminus \L}$ is a Riesz sequence for
$L^2(\T\setminus S)$.
 Moreover, in this case the sampling bound
in (\ref{frame}) and the Riesz sequence bound in (\ref{rs}) are
equivalent, $A/2\leq B\leq 2A$.
\end{clame}
Indeed, assume (\ref{frame}) holds.  
Let $P(t)=\sum_{\gamma\in\Z\setminus \L} a_{\gamma}e^{i\gamma t}$. Then
\[\|P\|^2_{L^2(S)}\leq
 \frac{1}{B}\sum_{\L}|\langle P,e^{i\l t}\rangle_{L^2(
S)}|^2=
\]
\[
\frac{1}{B}\sum_{\L}|\langle P,e^{i\l t}\rangle_{L^2
(\T\setminus S)}|^2\leq \frac{1}{B}\|P\|^2_{L^2(
\T\setminus S)},
\]
Since $B\leq 1$, this proves (\ref{rs}) with $A\geq B/2$,
$\Gamma=\Z\setminus \L$ and $\Omega=\T\setminus S$. The proof of
the opposite direction follows in much the same way.

\subsection{Proof of Theorem 2}
It suffices to show that Theorem 2 holds for open sets. Fix $d>0$
and let $S$ be an open set such that $(1+d)|S|<1$. Below we
denote by $C(d)$  positive constants, which depend
only on $d$.

Let $Q_j$ be disjoint intervals such that $S=\cup_{j\in \N} Q_j$.
For every $m\in \N$ denote $S_m:= \cup_{j<m}Q_j$. Using Theorem 1,
find for every $m$ a set $\L_m\subset \Z$, such that
$D(\L_m)<(1+d)|S_m|$ and $\L_m$ is a sampling sequence for
$PW_{S_m}$ with sampling bound  larger  then $C(d)|S_m|$.

For every $m$ choose a finite set $\Gamma_m$ and an interval
$I_m\subset \Z$ so that:
\[
\Gamma_m\subset (\Z\setminus \L_m)\cap I_m,
\]
and
\[
\frac{|\Gamma_m|}{|I_m|}>1-(1+d)|S_m|.
\]
By Claim 1 the set of exponentials $\{e^{i\gamma
t}\}_{\gamma\in\Gamma_m}$ is a Riesz sequence in $L^2(\T\setminus
S_m)$ with Riesz sequence bound larger then $C(d)|S_m|$.

Let $\textit{H}$ be the minimal translation invariant family of
finite sets in $\Z$, which contains all finite subsets of every
$\Gamma_m$. For this $\textit{H}$ take $\Gamma\subset\Z$ according
to Theorem B. Then
$$D^{\sharp}(\Gamma)\geq 1-(1+d){|S|}.$$
Take an arbitrary integer $M$. Then $\Gamma^M:=\Gamma\cap [-M,M]$
belongs to $\textit{H}$. So we can fix $m=m(M)$ such that
$\Gamma^M$ is a subset of some shifted $\Gamma_m$. Clearly $m$ tends
to infinity as $M$ does. So $\{e^{i\gamma t}\}_{\gamma \in
\Gamma^M}$ is a (finite) Riesz sequence in $L^2(\T\setminus S_m)$
with bound bigger then $C(d)|S_m|$. It follows $\{e^{i\gamma
t}\}_{\gamma \in \Gamma}$ is a Riesz in $L^2(\T\setminus S)$ with
bound bigger then $C(d)|S|$. Using claim 1 again we conclude that
$\L:=\Z\setminus\Gamma$ satisfies the required conditions.

\section{Upper bounds}

\subsection{Good Bessel sequences} In this section we discuss the possible upper bounds in
a sampling process. A sequence $\L\subset \Z$ is called a Bessel
sequence for $PW_S$ if there exists $B>0$ such that
\[
 \sum_{\l\in\L} |f(\l)|^2\leq B\|f\|^2
\:\:\:\:\:\:\:\:\:\:\:\:\forall f\in PW_S.
\]

Here we show that a theorem of Lunin \cite{Lu} implies the
existence of good Bessel sequences: Sequences which are not too
sparse on one side, but provide good Bessel bounds on the other
side. Again, we consider spectra of small measure. We keep in mind
the canonical example from Section 1 and show that the
inequalities in Theorems 1 and 2 can be reversed.

\begin{theorem}
 Given a compact set $S$, one can find
a sequence $\L$ of uniform density $D(\L)> |S|$ such that
\[
\sum_{\l\in\L} |f(\l)|^2\leq C|S|\|f\|^2 \:\:\:\:\:\:\:\:\:\:\:\:\forall
f\in PW_S,
\]
where $C$ is an absolute positive constant.
\end{theorem}
Moreover, in this case the sequence $\L$ can be constructed as a
finite union of arithmetic progressions.
\begin{theorem}
 Given a set $S$ of positive measure,
one can find a sequence $\L$ of density $D^{\sharp}(\L)> |S|$ such
that
\[
\sum |f(\l)|^2\leq C|S|\|f\|^2 \:\:\:\:\:\:\:\:\:\:\:\:\forall
f\in PW_S,
\]
where $C$ is an absolute positive constant.
\end{theorem}

 \subsection{Lunin's theorem} The type of problems studied in
 theorem A goes back to
Kashin, \cite{K1}, \cite{K2}, who in connection with certain
questions in the theory of orthogonal series proved that for every
$m\times n$ matrix $A$ of norm $1$ with $n/m<\rho(\eps)$, one can
find a subset $J\subset \{1,...,m\}$, of size $|J|=n$, so that the
sub-matrix $A(J)$ has norm less then $\eps$. Moreover, some
estimate of $\rho(\eps)$ was established. A sharp constant
$\rho(\eps)$ was found by Lunin \cite{Lu}:

\begin{theoremc}
Let $A$ be an $m\times n$ with $\|A\|=1$. Then there exists a
subset $J\subset \{1,...,m\}$ for which $|J|= n$ and
\[
\|A(J)w\|_{l^n_2}^2\leq
C\frac{n}{m}\|w\|^2\:\:\:\:\:\:\:\:\:\:\:\:\forall w\in l_2^n,
\]
where $C$ is an absolute positive constant.
\end{theoremc}

Observe that Theorem C can be proved also using the technique
developed in~\cite{S1}. 

One can \textit{prove Theorem 3}  in the same way as Theorem 1, the only difference is that one uses  Theorem C instead of Corollary 1.

The \textit{proof of Theorem 4} is pretty similar to the proof of
Theorem 2, in which one has to use Theorem 3 instead of Theorem 1.

\section{Additional remarks and open problems}

\subsection{Restricted Invertibility Theorem}
The Restricted Invertibility Theorem (RIT) was proved by Bourgain
and Tzafriri in \cite{BT} using a probabilistic approach, see also
\cite{BT2}, \cite{KT}, and \cite{V}. Recently Spielman and
Srivastava, \cite{S2}, used the technique developed in \cite{S1}
to give a short linear algebraic proof of the theorem. We
formulate here the RIT as in \cite{S2}.

\begin{theoremd} Fix $0<d<1$. For any operator $T:l_2^m\rightarrow l_2^m$ which satisfies
$\|Te_i\|=1$, where $\{e_i\}$ is the canonical basis for $l_2^m$,
one can find a subset $J\subset \{1,...,m\}$ of size
$|J|>(1-d){m}/{\|T\|^2}$ which satisfies
\[
C(d)\sum_{i\in J}|a_i|^2\leq \|\sum_{i\in J} a_iTe_i\|^2,
\]
where $C(d)=(1-\sqrt{1-d})^2$.
\end{theoremd}

In \cite{BT} Bourgain and Tzafriri applied the RIT to obtain a
result regarding good Riesz sequences for general measurable
sets. For an open set, this result may be improved by choosing
the frequencies of the corresponding Riesz sequence to be
uniformly distributed in $\Z$. To this end, we need the following
corollary of the RIT.

\begin{coro}
Fix $0<d<1$. Let $A$ be an $m\times n$ matrix which is a
sub-matrix of some $m\times m$ orthonormal matrix and such that
all of it's rows have equal $l^2$ norm. Then there exists a subset
$J\subset \{1,...,m\}$ for which $|J|\geq (1-d) n$ and
\[
C(d)\frac{n}{m}\|w\|^2\leq \|A(J)^{\mathrm{T}}
w\|_{l^n_2}^2\:\:\:\:\:\:\:\:\:\:\:\:\forall w\in l_2(J),
\]
where $C(d)$ is the constant given in Theorem D.
\end{coro}

Indeed, without loss of generality we may assume that the columns
of $A$ are the first $n$ columns of some $m\times m$ orthonormal
matrix. Let $P$ be the orthonormal projection on $l^2_m$ which is
defined by replacing all but the first $n$ coordinates of a vector
by zero. Then $T=\sqrt{m/n}P$ satisfies the conditions of Theorem
D, corresponding to the basis given by the rows of the orthonormal
matrix, and $\|T\|^2=m/n$. The corollary follows.

\begin{theorem}
Fix $0<d<1$. Given an open set $\Omega\subset[0,2\pi]$, there
exists $\Gamma\subset\Z$ of uniform density $D(\Gamma)>
(1-d){|\Omega|}$ such that
\[
C(d)|\Omega|\sum_{\gamma\in\Gamma}|a_{\gamma}|^2\leq \|\sum_{\gamma\in\Gamma} a_{\gamma}e^{i\gamma
t}\|_{L^2(\Omega)}^2,
\]
where $C(d)$ is the constant given in Theorem D.
\end{theorem}

This theorem can be proved using Corollary 2 in much the same way
that Theorem 1 is obtained from Corollary 1. We omit the proof.

\noindent\textbf{Remark 1.} Using Claim 1, one can reformulate this
result in terms of sampling. In particular, this implies the
existence of good sampling sequences in the sense of the
definition in Section~1, for sets whose measure is not too
small. Moreover, by an appropriate choice of the constant $d$ in
Theorem 5, one can  get in this way a weaker version of Theorem 1,
with $|S|$ replaced by $|S|^2$ in the sampling bound. However, it
is not clear whether Theorem 1 itself can be deduced this way. At
least, the value of the constant $C(d)$ in Theorem D does not
imply this.

\subsection{Open problems}
\begin{itemize}
\item[1.] It seems to be an important problem, whether Theorems 1 and 3 could be combined: Given a compact $S$ of positive measure, is there a good
exponential frame, that is a sequence $\L\subset\Z$ such that the
two sided inequality
\[
A|S|\|f\|^2\leq \sum | f(\l)|^2\leq
B|S|\|f\|^2\:\:\:\:\:\:\:\:\:\:\:\:\forall f\in PW_S
\]
is fulfilled with absolute constants $A$ and $B$? It is mentioned
in \cite{S1} that the finite dimensional result of this type would
imply the Kadison-Singer conjecture.

\item[2.] Another related problem: Given a compact $\Omega\subset[-\pi,\pi]$ of positive measure, is there a sequence\ $\Gamma\subset \Z$ without arbitrarily long gaps, such that
$\{e^{i\gamma t}\}_{\gamma\in\Gamma}$ is a Riesz sequence in
$L^2(\Omega)$? See \cite{LAW}, where it is proved that this
question is equivalent to the Kadison-Singer conjecture for
exponential frames.

\end{itemize}



\address{School of Mathematical Sciences, Tel Aviv University, Ramat-Aviv 69978, Israel
} \email{nitzansi@post.tau.ac.il}

\address{School of Mathematical Sciences, Tel Aviv University, Ramat-Aviv 69978, Israel
} \email{olevskii@math.tau.ac.il}

\address{University  of Stavanger,  Stavanger 4036, Norway}\\
\email{alexander.ulanovskii@uis.no}


\end{document}